\documentclass[11pt]{article}
\usepackage[a4paper, total={6in, 8in}]{geometry}
\usepackage{thmtools}
\usepackage{thm-restate}
\usepackage{authblk}
\usepackage{amsmath}
\usepackage{amsthm}
\usepackage[OT1]{fontenc}
\numberwithin{equation}{section}
\usepackage{amsfonts}
\usepackage{amssymb}
\usepackage{graphicx}
\usepackage{esint}
\usepackage{hyperref}
\usepackage{bm}
\usepackage{comment}
\hypersetup{linktocpage,
colorlinks=true,
linktoc=all,
linkcolor=blue,
}

\makeatletter

\newtheorem{theorem}{Theorem}[section]

\newtheorem{corollary}{Corollary}[section]
\newtheorem{definition}{Definition}[section]

\newtheorem{lemma}{Lemma}[section]

\usepackage[color=magenta]{todonotes}

 \newcommand{\D}{{\mathcal{D}}}

\newcommand{\Ra}{{\rm{Ra}}}
\newcommand{\Nu}{{\rm{Nu}}}
\newcommand{\Pe}{{\rm{Pe}}}
\newcommand{\ub}{{\bm{u}}}
\newcommand{\xx}{{\bm{x}}}
\newcommand{\yy}{{\bm{y}}}

\newcommand{\vv}{{\bm{v}}}

\title{Wall-to-wall heat transport by stochastic velocity fields}

\author[1]{Franco Flandoli}
\author[2]{Theresa Lange}
\author[3]{Camilla Nobili}
\affil[1]{Scuola Normale Superiore di Pisa, Piazza dei Cavalieri 7, 56126 Pisa, Italy}
\affil[2]{Max-Planck Institute for Mathematics in the Sciences, Inselstr. 22-26, 04103 Leipzig, Germany}
\affil[3]{School of Mathematics and Physics, University of Surrey, Guildford, Surrey GU2 7HX, UK}
\affil[ ]{E-mail: franco.flandoli@sns.it, theresa.lange@mis.mpg.de, c.nobili@surrey.ac.uk}
\begin{document}

\maketitle

\begin{abstract}
Inspired by the wall-to-wall optimal transport problem of Hassanzadeh,
Chini and Doering \cite{HCD2014}, we study heat transport between two isothermal plates
by an incompressible velocity field modelled as a transport noise: white
in time, coloured in space, and acting on the temperature in the
Stratonovich sense. The temperature is treated as a passive scalar, and
no momentum equation is imposed on the velocity. Passing to the It\^o
formulation, the mean temperature solves a deterministic diffusion
equation in which the noise induces an additional, spatially
inhomogeneous diffusion, and the Nusselt number $\Nu$ is expressed directly in
terms of the vertical covariance of the velocity field. We prove a
general upper bound on the $\Nu$ in terms of this covariance
near the walls, valid for every admissible noise. The no-slip boundary conditions
and incompressibility force the vertical covariance to decay
polynomially at the plates. When the intensity of the noise is measured
by a covariance budget $\Pe^2$ and its correlation structure does not
depend on $\Pe$, this decay yields $\Nu\lesssim\Pe^{2/3}$, and
$\Nu\lesssim\Pe^{1/2}$ when vertical and horizontal velocity components
are uncorrelated. We will show that this scaling is in fact sharp. Deterministic flows
under an energy budget can achieve $\Nu\sim\Pe$, while the exponent $2/3$
is the optimal one under an enstrophy budget \cite{doering2019optimal}. Transport noise thus
reaches this exponent from a budget of energy type, the derivative
control entering only through the behaviour of the noise at the walls,
and the mechanism is the diffusive structure of the It\^o-Stratonovich corrector.
The present setting in which the velocity is prescribed, is a first
step towards stochastic convection in which the noise is coupled to the
temperature through buoyancy.
\end{abstract}

\newpage
\section{Introduction}
\noindent
This manuscript investigates upper bounds on the heat transport between two
isothermal plates when the advecting velocity field is modelled as a
stochastic process, white in time and coloured in space, entering the
advection--diffusion equation as a transport noise in Stratonovich form.

We adopt a purely kinematic, passive-scalar point of view, in the spirit of the wall-to-wall optimal transport problem introduced by Hassanzadeh, Chini and Doering \cite{HCD2014} and later analysed rigorously by Tobasco and Doering \cite{TD2017, doering2019optimal}. The velocity field is not required to satisfy any momentum equation, and we ask how much heat can be transported between the plates by incompressible
stochastic velocity fields, given a \textit{budget on the covariance of the noise}, in the regime of large P\'eclet number $\Pe\gg1$ (see below).

In the box $[0,hL]^{d-1}\times[0,h]$, $d\geq 3$, we consider the
temperature $\Theta$ advected by an incompressible velocity field $\vv$,
\begin{equation*}
\partial_t\Theta + (\vv\cdot\nabla)\Theta=\kappa\Delta\Theta,\qquad
\nabla\cdot \vv=0,
\end{equation*}
periodic in $\xx' \in [0,hL]^{d-1}$, with $\Theta=\Theta_{\rm b}$ at $x_d=0$,
$\Theta=\Theta_{\rm t}$ at $x_d=h$, and no-slip boundary conditions $\vv=0$ at
$x_d\in\{0,h\}$. Rescaling lengths by $h$, time by $h^2/\kappa$, and
setting $T=(\Theta-\Theta_{\rm t})/(\Theta_{\rm b}-\Theta_{\rm t})$,
$\ub=h\vv/\kappa$, we obtain the dimensionless system on $\D:= [0,L]^{d-1}\times[0,1]$
\begin{align}
	\partial_t T + (\ub \cdot \nabla) T &= \Delta T,\label{eq:T}\\
	\nabla \cdot \ub &=0,\label{eq:divU}\\
	T(t,\xx',0) &= 1, \quad T(t,\xx',1) = 0,\quad \forall t \geq 0, \xx' \in [0,L]^{d-1},\label{eq:Tbc}\\
	\ub(t,\xx',0) &= 0, \quad \ub(t,\xx',1) = 0, \quad \forall t \geq 0, \xx' \in [0,L]^{d-1}.\label{eq:Ubc}
\end{align} 
Throughout, $\vv$ denotes the dimensional and
$\ub$ the dimensionless velocity.

\medskip

For a deterministic velocity field, the strength of advection can be measured, as in \cite{HCD2014}, by the P\'eclet number
\begin{equation}\label{Peclet}
\Pe^2=\overline{\int_0^1\langle|\ub|^2\rangle'\,{\rm d}x_d}\,,
\end{equation} 
and the heat transport by the Nusselt number
\begin{equation*}
\Nu=\overline{\int_0^1\langle(\ub T-\nabla T)\cdot \bm{e}_d\rangle'\,{\rm d}x_d}
 =-\overline{\langle\partial_dT\rangle'}\big|_{x_d=0},
\end{equation*}
where the second identity follows from the fact that the horizontally
averaged vertical flux is independent of $x_d$ (the notation for the
averages is introduced at the end of this section).

In \cite{HCD2014} the authors seek the steady, two-dimensional,
divergence-free velocity fields that maximise $\Nu$ under this energy
budget, subject only to the impermeability condition
$u_d|_{x_d=0,1}=0$. A direct application of the Cauchy--Schwarz
inequality and the maximum principle gives the rigorous bound
$\Nu\leq 1+\Pe/2$, valid for any such flow. To identify the optimal flows,
the authors derive the Euler--Lagrange equations of the constrained
maximisation problem and follow their solutions by numerical
continuation, starting from the exact solution available for small
$\Pe$, where the optimisers are arrays of square convection cells. For
large $\Pe$, guided by the numerics, they construct asymptotic solutions
by matched boundary-layer analysis. For a fixed aspect ratio $\Gamma$ of
the cells, the maximal transport scales as
$\Nu\sim\Gamma^{-2/3}\Pe^{2/3}$, attained by flows with a well-mixed
interior and strong tangential jets along the walls. Allowing the cells
to narrow as $\Gamma_{\rm opt}\sim \Pe^{-1/2}$ yields
$\Nu_{\rm MAX}\sim 0.1875\,\Pe$. The optimal flows thus saturate the rigorous bound up to a constant
factor. These results rely on asymptotic analysis and numerical evidence
rather than on a proof of optimality.

\smallskip
Under an enstrophy budget $\overline{\int_0^1\langle|\nabla \ub|^2\rangle\,{\rm d}x_d}=\Pe^2$,
the numerical optimum of \cite{HCD2014} for stress-free walls scales
approximately as $\Pe^{0.58}$.
The enstrophy-constrained problem was later revisited by Tobasco and
Doering \cite{TD2017, doering2019optimal}, who proved the rigorous upper bound
$\Nu\lesssim \Pe^{2/3}$, independently of the space dimension and of the
velocity boundary conditions, and constructed branching flows achieving
this scaling up to logarithmic corrections, thereby establishing it as
optimal.

\medskip

In this paper we replace the deterministic velocity field by a stochastic one,
following a scale-separation ansatz, with the goal of understanding the sensitivity of such scaling laws with respect to a transport noise. A turbulent fluid exhibits a wide
spectrum of coexisting spatial scales: a large-scale mean flow,
accompanied by smaller-scale eddies whose interplay determines the
mixing of transported quantities such as the temperature. In practice,
small scales can only be resolved down to the resolution of the
measuring device or of the numerical scheme; below this threshold only
their statistical features are accessible. A common paradigm in
stochastic fluid dynamics is to reintroduce these unresolved scales into
the model as a noise whose spatial structure encodes their statistics,
while their temporal decorrelation is idealised as instantaneous. This
is the viewpoint of Kraichnan's model of passive scalar advection
\cite{Kraichnan} and of space scale separation in the more recent theory of transport noise (for a recent overview refer to e.g. \cite{FL2023}). Let us point out here that this separation appears to only be reasonable in three or higher dimensions (in two dimensions, the presence of an inverse cascade potentially leads to larger structures created from smaller ones which renders the separation ansatz no longer suitable). In natural convection at high intensity, the
small-scale velocity fluctuations are generated by small, closely packed
parcels of hot and cold fluid, therefore a statistical description of the
velocity field appears to be natural. In the present work, however, we do not model
this feedback: the stochastic velocity field is prescribed, and the temperature is
treated as a passive scalar. Coupling the noise to the temperature is
the subject of our subsequent work in \cite{FLN2026b,FLN2026c}.

Accordingly, we take the dimensionless velocity to be
\begin{equation}\label{eq:stochU}
\ub=\sqrt{\mathcal{Q}_{\ub}}\,\dot W,
\end{equation}
i.e.\ at each position the velocity is driven by white noise $\dot W$
(on a probability space $(\Omega,\mathcal F,\mathbb P)$), whose spatial interactions at different points are determined by the covariance operator
$\mathcal{Q}_{\ub}$. We assume that $\mathcal{Q}_{\ub}$ has a matrix-valued kernel $Q_{\ub}(\xx,\yy)=[q_{i,j}(\xx,\yy)]_{i,j=1,\dots,d}$, chosen such that $\ub$ is divergence-free and satisfies the no-slip boundary condition \eqref{eq:Ubc} in the distributional sense. Heuristically,
$\mathbb E[u_i(t,\xx)u_j(s,\yy)]=q_{i,j}(\xx,\yy)\,\delta(t-s)$, so that
$q_{i,j}$ describes only the spatial structure of the velocity; the
distributional-in-time nature of $\dot W$ does not allow this identity
to be written rigorously, and the precise construction is given in
Section \ref{sec:stochFramework}. We write
$\widetilde Q_{\ub}(\xx)=[\widetilde q_{i,j}(\xx)]_{i,j=1,\dots,d}:=Q_{\ub}(\xx,\xx)$. Note that the no-slip boundary conditions force $\widetilde Q_{\ub}=0$ on $\{x_d=0\}\cup\{x_d=1\}$.

With this choice, \eqref{eq:T} becomes a stochastic advection--diffusion
equation, which we interpret in the Stratonovich sense. Transferring to the corresponding It\^{o} formulation enables us to take expectations and we find that the mean temperature $\mathbb E[T]$ solves the deterministic
parabolic problem
\begin{equation}\label{eq:meanT}
\partial_t\mathbb E[T]
=\nabla\cdot\left(\left({\rm Id}+\tfrac12\widetilde Q_{\ub}\right)\nabla\mathbb E[T]\right)
\end{equation}
with boundary conditions \eqref{eq:Tbc}. The additional term $\frac{1}{2}\nabla \cdot (\widetilde{Q}_\ub \nabla \mathbb{E}[T])$ stems from the It\^{o}-Stratonovich correction and encodes the additional (spatially inhomogeneous) diffusion induced by the noise, so that the matrix
$\mathrm{Id}+\frac12\widetilde Q_{\ub}$ plays the role of an
\emph{effective diffusivity}.
This leads us to define the
Nusselt number through the mean temperature,
\begin{equation}\label{def:stochasticNu}
\Nu:=-\overline{\langle\partial_d\mathbb E[T]\rangle'}\big|_{x_d=0},
\end{equation}
and we will show in Section \ref{sec:thmGeneralBound} that
\begin{equation}\label{eq:NuCovariance}
\Nu=\int_0^1\overline{\left\langle-\partial_d\mathbb E[T]
-\frac{1}{2}\sum_{j=1}^d\widetilde q_{d,j}\,\partial_j\mathbb E[T]\right\rangle'}\,{\rm d}x_d.
\end{equation}
Identity \eqref{eq:NuCovariance} establishes a transparent link between
the heat transport and the vertical covariance $\widetilde q_{d,j}$, which
plays the role of the convective flux $\overline{\langle u_dT\rangle'}$
in the deterministic setting. While the It\^o--Stratonovich corrector of
transport noise has been exploited in other problems in fluid dynamics (e.g. refer to \cite{galeati2020convergence,flandoli2021high,flandoli2021delayed,flandoli2022eddy,agresti2026absence}),
to the best of our knowledge this relation between the noise covariance
and the Nusselt number has not previously been written down in the
context of stochastic convection.

Nusselt number bounds under stochastic forcing have also been investigated in  Földes, Glatt-Holtz, Richards
and Whitehead \cite{FGHRW2016} (for related works see also \cite{FGHRT2015}, \cite{FGHR2019}), but in a substantially different setting.
In particular, in  \cite{FGHRW2016}, the velocity is not prescribed but solves the Boussinesq
equations, and the randomness enters as an additive forcing of the
temperature (and possibly the momentum) equation. The Nusselt number is
defined as an average with respect to the unique invariant measure of
the system, whose existence and ergodicity are established for large
Prandtl number $\Pr$. The background-field method then yields
$\Nu\lesssim_{\D}\Ra^{1/2}$, where $\Ra$ is the Rayleigh number.\footnote{In \cite{FGHRW2016}, the authors adopt the following non-standard notation: $\Ra$ is a Rayleigh number built on the temperature
scale set by the noise intensity, $\widetilde{\Ra}$ is the ratio of the
imposed temperature difference to this scale, and the classical Rayleigh
number is the product $\Ra\,\widetilde{\Ra}$ (see the appendix of
\cite{FGHRW2016}). Their bound $\Nu\lesssim(\Ra\,\widetilde{\Ra})^{1/2}$ thus
recovers the classical $\Ra^{1/2}$ scaling, up to an additional It\^o
correction term.} In our setting of transport noise, such machinery is not needed: by \eqref{eq:meanT} the mean temperature obeys a closed deterministic
equation, thus $\Nu$ is well defined without appealing to invariant
measures.

\vspace{0.5cm}
To derive such bounds without enforcing a specific momentum equation (such as
the Stokes or Navier--Stokes equations), we work with an admissible class
of velocity fields characterised solely by incompressibility, the no-slip
boundary conditions and spatial regularity of the covariance. Our first result is a general upper bound on $\Nu$ in terms of the vertical covariance in a boundary layer of arbitrary thickness $\delta$.

\begin{theorem}\label{thm:NuGeneralBound}
Let $\ub$ be a stochastic velocity field of the form \eqref{eq:stochU},
divergence-free \eqref{eq:divU} and satisfying the no-slip boundary
conditions \eqref{eq:Ubc}. Then, for any $\delta\in(0,1]$, the Nusselt
number \eqref{def:stochasticNu} satisfies
\begin{equation}\label{eq:genUpperBound}
\Nu\lesssim_d \frac1\delta+\frac1{\delta^2}\sum_{j=1}^d\int_0^\delta
\big\langle|\widetilde q_{d,j}|^2\big\rangle'\,{\rm d}x_d ,
\end{equation}
where $\lesssim_d$ denotes an inequality up to a universal constant $C_d$
depending only on the spatial dimension $d$.
\end{theorem}
It remains to define the P\'eclet number in the stochastic case. The following definition will be adopted throughout the paper, replacing \eqref{Peclet}.
\begin{definition}
Let $\ub$ be a stochastic velocity field of the form \eqref{eq:stochU}.
We define the associated P\'eclet number as
\begin{equation}\label{stochastic_budget}
\Pe^{2}=\int_{0}^{1}\left\langle {\rm{tr}}\left(  \widetilde{Q}_\ub\right)
\right\rangle'{\rm d}x_d
\end{equation}
\end{definition}

Since $\ub$ is white in time, its kinetic energy is infinite, and the
fixed-energy constraint of \cite{HCD2014} has no direct analogue. The
natural quantity is instead the covariance at coincident points, which,
by \eqref{eq:meanT}, acts on the mean temperature as an additional
diffusion coefficient. To make this precise, we return to dimensional
variables and write the noise as $\sqrt{\mathcal{Q}_{\vv}}\,\dot{W}$. Since
$W$ has the units of square-root of time, the covariance kernel at
coincident points $\widetilde Q_{\vv}$ has the units of a diffusivity, and
its average over the layer,
\begin{equation}\label{budget-dimensional}
K:=\frac1h\int_0^h\left\langle\mathrm{tr}\big(\widetilde Q_{\vv}\big)\right\rangle'\,{\rm d}x_d ,
\end{equation}
measures the mean diffusion induced by the noise. We compare it with the
molecular diffusivity $\kappa$, setting
\begin{equation*}
\Pe^2:=\frac{K}{\kappa}=\frac{h^2/\kappa}{h^2/K}.
\end{equation*}
That is, the ratio between the time needed by molecular diffusion and
the time needed by the noise-induced diffusion to mix the layer, just
as, in the deterministic setting, $\Pe$ is the ratio between the
diffusive and the advective time. A rescaling of the equation, carried out in Appendix \ref{app:rescaling},
gives $\widetilde Q_{\ub}(\hat{\boldsymbol x})=\widetilde Q_{\vv}(h\hat{\boldsymbol x})/\kappa$:
the noise does not rescale like a velocity, since Brownian increments
 scale like the square-root of time. 
\medskip

Our goal is to quantify how fast the vertical covariance decays near the
walls, and here the boundary conditions play the decisive role. Writing
$\sqrt{\mathcal{Q}_{\ub}}\,\dot W=\sum_k\lambda_k\,e_k\,\dot W^k$, so that
$\widetilde q_{d,j}=\sum_k\lambda_k^2\,e_k^d\,e_k^j$ (for a precise formulation, c.f. Section \ref{sec:stochFramework}), the no-slip condition
imposes $e_k=0$ at $x_d=0$, and incompressibility then gives
$\partial_de_k^d=-\sum_{j<d}\partial_je_k^j=0$ there. Hence the vertical
components vanish quadratically, $e_k^d=O(x_d^2)$, while the horizontal
ones vanish linearly, so that $\widetilde q_{d,j}=O(x_d^3)$ for $j<d$ and
$\widetilde q_{d,d}=O(x_d^4)$. We stress that impermeability alone would only
give $e_k^d=O(x_d)$, thus the no-slip condition is essential.

The constants in these expansions involve derivatives of the $e_k$ at
the walls, on which the budget \eqref{stochastic_budget} gives no
control. To obtain a scaling law in $\Pe$, we let the intensity of the
noise grow while keeping its spatial structure fixed: we assume
\begin{equation}\label{eq:scalingHyp}
\lambda_k=\Pe\,\mu_k,\qquad
\sum_k\mu_k^2\int_0^1\big\langle|e_k|^2\big\rangle'\,{\rm d}x_d=1,
\end{equation}
where the divergence-free fields $e_k$ vanish at the walls, and the
coefficients $\mu_k$ do not depend on $\Pe$. 

Under \eqref{eq:scalingHyp}, the covariance kernel of the noise reads
\begin{equation}\label{eq:covarianceHyp}
Q_{\ub}(\boldsymbol x,\boldsymbol y)
=\Pe^2\sum_k\mu_k^2\,e_k(\boldsymbol x)\otimes e_k(\boldsymbol y).
\end{equation}
Taking the trace on the diagonal $\boldsymbol y=\boldsymbol x$, the
normalisation in \eqref{eq:scalingHyp} ensures that the budget
\eqref{stochastic_budget} is satisfied, as we can easily see from the identity
\[\mathrm{tr}\big(\widetilde Q_{\ub}\big)=\Pe^2\sum_k\mu_k^2|e_k|^2.\]
We stress that
\eqref{eq:scalingHyp} does not fix the velocity field: its realisations
are random combinations of the modes $e_k$, with coefficients
fluctuating in time. What is fixed is only its correlation structure,
i.e.\ the covariance \eqref{eq:covarianceHyp} up to the factor $\Pe^2$.
\begin{corollary}\label{cor:PeScaling}
Under \eqref{eq:scalingHyp}, for all $x_d\in[0,1]$,
\begin{align}
\left\langle|\widetilde q_{d,j}|^2\right\rangle'&\leq C\,\Pe^4x_d^6,
\qquad j=1,\dots,d-1, \label{eq:decay1}\\
\left\langle|\widetilde q_{d,d}|^2\right\rangle'&\leq C\,\Pe^4x_d^8,
\label{eq:decay2}
\end{align}
with $C$ depending only on $d$ and $\sum_k\mu_k^2\,\|e_k\|_{C^3}^2<\infty$. Consequently,
\[\Nu\lesssim\Pe^{2/3} \quad \mbox{ for }\quad \Pe\ge1.\]
If moreover
$\left\langle|\widetilde q_{d,j}|^2\right\rangle' = 0$ for $j<d$, then \[\Nu\lesssim\Pe^{1/2}\,,\] and this bound is sharp for $d=3$.
\end{corollary}

Let us comment on these results in light of the deterministic problem.
Under an energy budget \eqref{Peclet}, the optimal steady flows
of \cite{HCD2014} saturate the linear bound $\Nu\lesssim\Pe$ up to a
constant factor, as found by asymptotic analysis and numerics, by means
of convection cells that narrow as $\Pe$ grows and carry strong
tangential jets along the walls. Such flows have no fixed spatial
structure and are therefore excluded by hypothesis \eqref{eq:scalingHyp}.
The exponent $2/3$ in Corollary \ref{cor:PeScaling} instead coincides
with the rigorous bound obtained by Tobasco and Doering~\cite{TD2017} for
deterministic flows under an enstrophy budget, which is known to be
optimal up to logarithmic corrections. Our budget
\eqref{stochastic_budget} is of energy type, while the constants in
\eqref{eq:decay1}--\eqref{eq:decay2} involve derivatives of the
covariance at the walls; hypothesis \eqref{eq:scalingHyp} therefore
provides a control on the derivatives of the noise, akin to an enstrophy
budget, although only its behaviour near the walls enters the estimates.
The mechanism that renders this sufficient is the It\^o corrector: by
\eqref{eq:NuCovariance}, the stochastic convective flux
$\frac12\widetilde q_{d,j}\,\partial_j\mathbb E[T]$ is of diffusive form,
which allows it to be estimated directly against the gradient of the
mean temperature. The deterministic flux $u_dT$ has no such structure,
and our argument does not transfer directly to deterministic flows.

\medskip
Furthermore, in the proof of Corollary \ref{cor:PeScaling} we determine $e_1,e_2$ and $C_{\Pe}$ such that the three-dimensional velocity field
\[\ub= \sum_{k=1}^2 \lambda_k e_k \dot{W},\]
 with $\lambda_k \propto\Pe$ dictated by the budget
\eqref{stochastic_budget}
(see \eqref{eq:uSharp} for the exact form), achieves $\Nu\sim \Pe^{\frac 12}$. This shows that the exponent $\frac 12$ is sharp for $d=3$.
\medskip

The present work is the first step towards a more challenging problem,
in which the stochastic velocity field is no longer prescribed but is
coupled to the temperature through buoyancy, as in Rayleigh--B\'enard
convection driven by transport noise. In that setting the covariance of
the noise is itself determined by the dynamics, and the bounds derived
here serve as a kinematic benchmark against which the effect of the
feedback can be measured. This problem will be addressed in a
forthcoming work.

\paragraph{Notations.}

We will use the following notions:
\begin{itemize}
	\item long-time average
	\[\overline{(\cdot)} := \limsup_{t_0\rightarrow \infty}\frac{1}{t_0} \int_0^{t_0} (\cdot){\rm d}t;\]
	\item horizontal-space average
	\[\langle (\cdot)\rangle' := \frac{1}{L^{d-1}} \int_{[0,L]^{d-1}} (\cdot) {\rm d} \xx'.\]
\end{itemize}

\paragraph{Organization of the paper.}
A first crucial step of the paper is to rigorously set up the model of the velocity field as in \eqref{eq:stochU}, as done in Section \ref{sec:stochFramework}. This will allow us to properly define stochastic solutions of the system \eqref{eq:T}-\eqref{eq:Ubc}, and in particular solutions of \eqref{eq:meanT}. In Section \ref{sec:preliminaries} this set-up allows us to properly define the Nusselt number \eqref{def:stochasticNu} as well as related useful identities. The proofs of our main results Theorem \ref{thm:NuGeneralBound} and Corollary \ref{cor:PeScaling} are then to be found in Section \ref{sec:thmGeneralBound} and Section \ref{sec:boundPeclet}, respectively.
\section{Stochastic framework}\label{sec:stochFramework}
The formal object
\[\ub = \sqrt{\mathcal{Q}_\ub}\dot{W}\]
shall denote a cylindrical Brownian motion on a probability space $(\Omega, \mathcal{F}, \mathbb{P})$ with respect to a divergence-free complete orthonormal system (CONS), with spatial covariance operator $\mathcal{Q}_\ub$, i.e.
\begin{equation}\label{eq:cylinBM}
    {\rm d} \ub = \sqrt{\mathcal{Q}_\ub}(\xx){\rm d}W_t = \sum_{k \in \mathbb{N}} \lambda_k e_k(\xx) {\rm d}\beta_t^k
\end{equation}
where $(\lambda_k) \in \ell^2$, $(e_k)_k$ are a complete orthonormal system of $L_\sigma^2(\D) = \overline{\{f \in C_c^\infty: \nabla \cdot f =0\}}^{L^2(\D; \mathbb{R}^d)}$ with $e_k \in C(\overline{\D}; \mathbb{R}^d)$ and $ e_k(x_d=0,1) = 0$, while $\beta^k$ are independent one-dimensional standard Brownian motions. The kernel $Q_\ub$ of $\mathcal{Q}_\ub$, defined by
\[Q_\ub(\xx,\yy) = \sum_k \lambda_k^2 e_k(\xx) \otimes e_k(\yy),\]
therefore has the following properties:
\begin{itemize}
	\item symmetric, positive semi-definite;
	\item $Q_\ub \in L^2(\D \times \D)$;
	\item $\widetilde{Q}_\ub(\cdot) := Q_\ub(\cdot,\cdot) \in L^1(\D)$;
	\item $Q_\ub(\xx,\yy) = 0$ for $\xx \in \partial \D$ or $\yy \in \partial\D$.
\end{itemize}
Such a velocity field gives rise to the following stochastic advection-diffusion equation for the temperature field
\begin{equation*}
	{\rm d} T + \left(\left(\sqrt{\mathcal{Q}_\ub}\circ {\rm d}W_t\right)\cdot \nabla \right) T = \Delta T {\rm d}t
\end{equation*}
in Stratonovich form, or equivalently
\begin{equation}\label{eq:Ito}
	{\rm d} T + \left(\left(\sqrt{\mathcal{Q}_\ub} {\rm d}W_t\right)\cdot \nabla \right) T = {\rm div}\left(\left({\rm Id} +\frac{1}{2} \widetilde{Q}_\ub \right) \nabla T\right){\rm d}t
\end{equation}
in It\^{o} form (e.g. c.f. Section 3.4 in \cite{CG24}). We define corresponding solutions in the following way: first of all let $\widetilde{Q}_\ub(\xx) = [\widetilde{q}_{i,j}(\xx)]_{i,j=1,\dots,d}$ and assume that
\begin{equation}\label{eq:regTildeQ}
    \widetilde{Q}_\ub \in W^{1,\infty}(\D).
\end{equation}
For the following spaces
\begin{align*}
        H &:= H_0^1(\D), \text{ with inner product } \langle \cdot, \cdot \rangle_H\\
        V &:= H^2(\D) \cap H_0^1(\D),
    \end{align*}
one can show that the bilinear form
\[B(u,v) := \int_\D \left \langle \left({\rm Id} + \frac{1}{2}\widetilde{Q}_\ub(\xx)\right)(\nabla u)(\xx), (\nabla v)(\xx) \right\rangle {\rm d}\xx, \quad u,v \in V,\]
is coercive and bounded on $V$, hence induces a self-adjoint operator $\mathcal{A}$ with domain
\[D(\mathcal{A}) = \left\{ u \in V: \exists h \in H, B(u,v) = - \langle h, v \rangle_H \quad \forall v \in V \right\} \subset V\]
such that
\[\left \langle \mathcal{A} u, u \right \rangle_H = - B(u,u) \leq 0 \quad \forall u \in D(\mathcal{A}).\]
Hence $- \mathcal{A}$ generates a strongly continuous semi-group $e^{t\mathcal{A}}, t \geq 0$, on $H$.
\begin{definition}[Definition of probabilistic solution]\label{def:solT}
	Let $\tau \geq 0$ and let $\mathfrak{T}= \mathfrak{T}(x_d)$ be a smooth deterministic function satisfying $\mathfrak{T}(0)=1, \mathfrak{T}(1)=0$. Given $T_0$ such that $\mathcal{T}_0 := T_0 - \mathfrak{T} \in H$, we say that $T$ is a \emph{probabilistically strong, mild solution} to \eqref{eq:Ito} with initial condition $T_0$, if $\mathcal{T}:= T - \mathfrak{T}$ is progressively measurable with paths of class $L^\infty(0,\tau;H) \cap L^2(0,\tau;V)$, and satisfies
	\begin{equation}\label{eq:mildSol}
		\mathcal{T}(t) = e^{t\mathcal{A}}\mathcal{T}_0 - \int_0^t e^{(t-s)\mathcal{A}}\left(\sqrt{\mathcal{Q}_\ub} {\rm d} W_s\right)\cdot \nabla \mathcal{T}(s) - \int_0^t e^{(t-s)\mathcal{A}}\left(\sqrt{\mathcal{Q}_\ub} {\rm d} W_s\right)\cdot \mathfrak{T}'\bm{e}_{x_d}
	\end{equation}
	with homogeneous boundary conditions, as well as the energy bound
	\begin{equation}\label{eq:energyInequ}
		\mathbb{E}\left[\sup_{s \leq t}\|\mathcal{T}(s)\|^2_{L^2(\D)}\right] + \mathbb{E}\left[\int_0^t \left\|\nabla \mathcal{T}(s)\right\|_{L^2(\D)}^2{\rm d}s\right] < \infty.
	\end{equation}
\end{definition}
\noindent
The existence of such solutions is rather classical and we refer e.g. to Chapter 3 in \cite{FL2023}, and references therein. Furthermore we observe that the expected value of the It\^{o}-integrals in \eqref{eq:mildSol} vanishes due to square integrability coming from \eqref{eq:energyInequ} thus following from our assumptions on the vector fields $(e_k)_k$ and using the regularity assumption \eqref{eq:regTildeQ} we deduce that $E:= \mathbb{E}[T]$ is a strong solution to
	\begin{align}
		&\partial_t E =\nabla \cdot\left(\left({\rm Id} + \frac{1}{2}\widetilde{Q}_\ub\right)\nabla E\right),\label{eq:EvolExp}\\
        &E(t,\xx',0) =1, \quad E(t,\xx',1) = 0, \label{eq:EvolExpBC}
	\end{align}
in the sense that $\mathcal{E} = E - \mathfrak{T}$ satisfies
\[\mathcal{E}(t) = e^{t \mathcal{A}}\mathcal{E}(0)\]
with homogeneous boundary conditions, and
\[\mathcal{E} \in L^\infty(0,\tau;H)\cap L^2(0,\tau;V), \,\partial_t \mathcal{E} \in L^2(0,\tau;L^2(\D)).\]
\section{Preliminaries}\label{sec:preliminaries}
\noindent
\begin{lemma}
    The Nusselt number defined by
    \begin{equation}\label{eq:stochasticNu}
        \Nu := -\left.\overline{\left\langle \partial_d \mathbb{E}[T]\right\rangle' }\right|_{x_d =0}
    \end{equation}
    is well-defined.
\end{lemma}
\begin{proof}
    Denote $E := \mathbb{E}[T]$. Then we know that 
    \begin{equation}\label{eq:reg1}
        \left({\rm Id} + \frac{1}{2}\widetilde{Q}_\ub\right)\nabla E \in L^2(0,\tau; H^1(\D)).
    \end{equation}
    Taking the horizontal average and using periodicity we obtain
    \begin{equation}\label{eq:horAvE}
        \partial_t \langle E\rangle' = \partial_d \left\langle\left[\left({\rm Id} + \frac{1}{2} \widetilde{Q}_\ub\right)\nabla E\right]^{(d)}\right\rangle'.
    \end{equation}
    From \eqref{eq:reg1} we obtain that
    \[\left \langle \left[\left({\rm Id} + \frac{1}{2}\widetilde{Q}_\ub\right)\nabla E\right]^{(d)} \right\rangle' \in L^2(0,\tau; H^1(0,1)),\]
    i.e.
    \[\int_0^\tau \left\|\left \langle \left[\left({\rm Id} + \frac{1}{2}\widetilde{Q}_\ub\right)\nabla E\right]^{(d)} \right\rangle'\right\|^2_{H^1}{\rm d} t < \infty.\]
    In particular this implies that
    \[\overline{\left \langle \left[\left({\rm Id} + \frac{1}{2}\widetilde{Q}_\ub\right)\nabla E\right]^{(d)} \right\rangle'}< \infty.\]
    Now using the Sobolev embedding $H^1(0,1) \hookrightarrow C(0,1)$ we are allowed to evaluate at $x_d =0$ and obtain
   \[\left.\overline{\left \langle \left[\left({\rm Id} + \frac{1}{2}\widetilde{Q}_\ub\right)\nabla E\right]^{(d)}\right \rangle '}\right|_{x_d =0} = \left.\overline{\left \langle \partial_dE\right \rangle '}\right|_{x_d =0} =: -\Nu\]
   by the boundary assumption on $\widetilde{Q}_\ub$. This concludes the proof.
\end{proof}
Next we point out that one can show that $E$ satisfies a maximum principle, i.e. if $\mathbb{E}[T_0] \in [0,1]$, then $\|E\|_{L^{\infty}_{t,\xx}} \leq 1$. Therefore we can show that the following identiy holds:
\begin{lemma}[Nusselt number identity]\label{lem:NuId}
	The Nusselt number as defined in \eqref{eq:stochasticNu} satisfies
	\begin{equation*}
		\Nu = \overline{\int_0^1 \left \langle \left| \left( {\rm Id} + \frac{1}{2}\widetilde{Q}_\ub\right)^{\frac{1}{2}}\nabla \mathbb{E}[T]\right|^2\right \rangle'{\rm d}x_d}.
	\end{equation*}
\end{lemma}
\begin{proof}
	Testing \eqref{eq:EvolExp} against $E$ and applying the spatial horizontal average, we obtain
    \begin{equation*}
        \frac{1}{2}\left \langle  \partial_t E^2 \right \rangle' = \left \langle \nabla \cdot\left(\left({\rm Id} +\frac{1}{2}\widetilde{Q}_\ub \right)\nabla E\right) E\right \rangle'.
    \end{equation*}
	It holds
	\begin{align*}
	    \left \langle \nabla \cdot\left(\left({\rm Id} +\frac{1}{2}\widetilde{Q}_\ub \right)\nabla E\right) E \right \rangle'
        &= \left \langle \nabla \cdot\left(\left({\rm Id} +\frac{1}{2}\widetilde{Q}_\ub \right)(\nabla E) E\right)\right \rangle' - \left \langle \left|\left({\rm Id} + \frac{1}{2}\widetilde{Q}_\ub\right)^{\frac{1}{2}}\nabla E\right|^2\right \rangle',
	\end{align*}
	hence
	\begin{align*}
	    \frac{1}{2}\left \langle \partial_t E^2 \right \rangle' &=  \partial_d \left \langle E \left[\left({\rm Id} + \frac{1}{2}\widetilde{Q}_\ub\right)\nabla E\right]^{(d)}\right \rangle' - \left \langle \left|\left({\rm Id} + \frac{1}{2}\widetilde{Q}_\ub\right)^{\frac{1}{2}}\nabla E\right|^2\right \rangle'.
	\end{align*} 
	Integrating in $x_d$ and time we obtain
	\begin{align*}
		0 &= \frac{1}{2} \int_0^t \int_0^1 \left \langle  \partial_t E^2 \right \rangle' {\rm d}x_d {\rm d}s - \int_0^t \int_0^1 \partial_d \left \langle E \left[\left({\rm Id} + \frac{1}{2}\widetilde{Q}_\ub\right)\nabla E\right]^{(d)}\right \rangle' {\rm d} x_d {\rm d}s\\
		&\quad + \int_0^t \int_0^1 \left \langle \left|\left({\rm Id} + \frac{1}{2}\widetilde{Q}_\ub\right)^{\frac{1}{2}}\nabla E\right|^2\right \rangle' {\rm d} x_d {\rm d} s\\
		&= \frac{1}{2} \int_0^t \int_0^1 \left \langle  \partial_t E^2 \right \rangle' {\rm d}x_d {\rm d}s + \int_0^t \left.\left \langle E \left[\left({\rm Id} + \frac{1}{2}\widetilde{Q}_\ub\right)\nabla E\right]^{(d)}\right \rangle'\right|_{x_d=0} {\rm d}s\\
		&\quad + \int_0^t \int_0^1 \left \langle \left|\left({\rm Id} + \frac{1}{2}\widetilde{Q}_\ub\right)^{\frac{1}{2}}\nabla E\right|^2\right \rangle' {\rm d} x_d {\rm d} s
	\end{align*}
	where we used the boundary conditions \eqref{eq:EvolExpBC}. By the maximum principle we obtain that $\int_0^t \left \langle  \partial_t E^2 \right \rangle' {\rm d}s$ is bounded. Further by the imposed boundary conditions on $E$ and $\widetilde{Q}_\ub$, respectively, we have
	\[\left.\left \langle E \left[\left({\rm Id} + \frac{1}{2}\widetilde{Q}_\ub\right)\nabla E\right]^{(d)}\right \rangle'\right|_{x_d=0} = \left.\left \langle \partial_d E\right \rangle'\right|_{x_d=0},\]
	hence we deduce the claim by definition of the Nusselt number in \eqref{eq:stochasticNu}.
\end{proof}
As this will become relevant later, let us finally discuss the following identity: recall \eqref{eq:horAvE} and consider the case $\widetilde{q}_{d,j} = \widetilde{q}_{d,j}(x_d)$, $j=1,\dots,d$. By periodicity in the horizontal variables (consequently $\langle \partial_j E\rangle'=0, j < d$), \eqref{eq:horAvE} reduces to the following closed system
\begin{align}
    \partial_t F &= \partial_3\left(\left(1+\frac{1}{2}\widetilde{q}_{d,d}(x_d)\right)\partial_d F\right),\label{eq:DiffEqu}\\
    F(t,\xx', 0) &= 1, \quad
    F(t,\xx', 1) = 0 \qquad \forall t \geq 0\label{eq:DiffEquBC},
\end{align}
where we may choose $F(0,\cdot) = F_0 \in L^2(0,1)$ with $F_0 \equiv 0$. In this case we can show the following identity:
\begin{lemma}\label{lem:NuReprZdep}
    In the case $\widetilde{q}_{d,j} = \widetilde{q}_{d,j}(x_d)$, $j=1,\dots,d$, the Nusselt number has the explicit representation
    \begin{equation}\label{eq:NuIdentity}
    \Nu = \left(\int_0^1 \frac{1}{1+ \frac{1}{2}\widetilde{q}_{d,d}(x)}{\rm d} x\right)^{-1}.
\end{equation}

\end{lemma}
\begin{proof}
    Denote $a(x) := 1+ \frac{1}{2}\widetilde{q}_{d,d}(x)$. The time-stationary solution to \eqref{eq:DiffEqu}-\eqref{eq:DiffEquBC} is
\begin{equation}
    f(x) = 1 - \frac{\int_0^x a(y)^{-1} {\rm d}y}{\int_0^1 a(y)^{-1} {\rm d}y}.
\end{equation}
In particular, it holds $f \in H^1(0,1).$ Then $\mathcal{F} := F - f$ solves the corresponding homogeneous system
\begin{align*}
    \partial_t \mathcal{F} &= \partial_d\left(a(x_d)\partial_d \mathcal{F}\right),\\
    \mathcal{F}(t,\xx',0) &= \mathcal{F}(t,\xx',1) = 0 \qquad \forall t \geq 0,\\
    \mathcal{F}(0,\cdot) &= -f.
\end{align*}
By classical theory, there exists a unique self-adjoint operator $\mathcal{A}^\mathcal{S}$ such that $-\mathcal{A}^\mathcal{S}$ generates a strongly continuous semi-group $(\mathcal{S}(t))_{t \geq 0}$ on $L^2(0,1)$ such that
\begin{equation}
    \mathcal{F}(t) = \mathcal{S}(t) f.
\end{equation}
Observe that for all $g\in L^2(0,1)$ it holds $\mathcal{S}(t)g \longrightarrow 0$ as $t \rightarrow \infty$. Further one can show that the operator
\[{\rm tr}: D\left(\mathcal{A}^\mathcal{S}\right) \rightarrow \mathbb{R}, g \mapsto \partial_x g(0)\]
is well-defined and bounded. Consequently it holds
\begin{equation}
    \lim_{t \rightarrow\infty}\frac{1}{t}\int_0^t \partial_d \left.\mathcal{F} \right|_{x_d = 0}{\rm d}s = 0.
\end{equation}
Indeed, by definition it holds
\[(\partial_d \mathcal{F})(s,0) = {\rm tr}(\mathcal{F}(s, \cdot)),\]
hence (w.l.o.g. $s>1$)
\begin{align*}
    \left|(\partial_d \mathcal{F})(s,0)\right| &\leq \left| {\rm tr}(\mathcal{F}(s, \cdot))\right|\\
    &\lesssim \|\mathcal{F}(s,\cdot)\|_{D\left(\mathcal{A}^\mathcal{S}\right)} \simeq \|\mathcal{S}(s)f\|_{D\left(\mathcal{A}^\mathcal{S}\right)} \simeq \|\mathcal{S}(1) \mathcal{S}(s-1)f\|_{D(\mathcal{A})}\\
    &\lesssim  \|\mathcal{S}(s-1)f\|_{L^2(0,1)} \longrightarrow 0, s \rightarrow \infty.
\end{align*}
Therefore we obtain the claimed Nusselt number identity
\begin{align*}
    \Nu = -  \left.\overline{\partial_d \langle E \rangle'} \right|_{x_d = 0}
    &= - \left.\overline{\partial_d F}\right|_{x_d = 0}\\
    &= - \left.(\partial_d f)\right|_{x_d=0} - \left.\overline{\left( \partial_d\mathcal{F}\right)}\right|_{x_d = 0}\\
    &=  \left(\int_0^1 a(x_d)^{-1} {\rm d} x_d\right)^{-1}.
\end{align*}
\end{proof}

\section{Proof of Theorem \ref{thm:NuGeneralBound}}\label{sec:thmGeneralBound}
Our main result makes use of the following exchangeability result (the proof will be provided at the end of this section):
\begin{lemma}\label{lem:ExOp2}
	There exists a constant $C$ independent of $x_d$ such that
	\begin{equation*}
		C = \overline{\partial_d \langle \mathbb{E}[T]\rangle' + \frac{1}{2} \left \langle \left[\widetilde{Q}_\ub \nabla \mathbb{E}[T]\right]^{(d)}\right\rangle'}.
	\end{equation*}
    Thus for any $\delta \in (0,1]$ it holds
	\begin{equation*}
		\frac{1}{\delta} \int_0^{\delta} \overline{\partial_d \langle \mathbb{E}[T] \rangle' + \frac{1}{2} \left\langle \left[\widetilde{Q}_\ub \nabla \mathbb{E}[T]\right]^{(d)}\right\rangle'} {\rm d}x_d =  \overline{\frac{1}{\delta} \int_0^{\delta} \partial_d \langle \mathbb{E}[T] \rangle' + \frac{1}{2} \left\langle \left[\widetilde{Q}_\ub \nabla \mathbb{E}[T]\right]^{(d)}\right\rangle' {\rm d} x_d}.
	\end{equation*}
\end{lemma}
\noindent
Recall that $\widetilde{Q}_\ub = [\widetilde{q}_{i,j}]_{i,j = 1,..., d}$, then with the previous result we can prove the main theorem, which we restate, for convenience of the readers:
\begin{theorem}
Let $\ub$ be a stochastic velocity field of the form \eqref{eq:stochU},
divergence-free \eqref{eq:divU} and satisfying the no-slip boundary
conditions \eqref{eq:Ubc}. Then, for any $\delta\in(0,1]$, the Nusselt
number \eqref{def:stochasticNu} satisfies
\begin{equation*}
\Nu\lesssim_d \frac1\delta+\frac1{\delta^2}\sum_{j=1}^d\int_0^\delta
\big\langle|\widetilde q_{d,j}|^2\big\rangle'\,{\rm d}x_d ,
\end{equation*}
where $\lesssim_d$ denotes an inequality up to a universal constant $C_d$
depending only on the spatial dimension $d$.
\end{theorem}
\begin{proof}
	For a solution $T$ as in Definition \ref{def:solT}, again denote $E := \mathbb{E}[T]$. From Lemma \ref{lem:ExOp2} we obtain a constant $C$ independent of $x_d$ such that
	\begin{equation*}
		C = \overline{\partial_d \langle E\rangle' + \frac{1}{2} \left \langle \left[\widetilde{Q}_\ub \nabla E\right]^{(d)}\right\rangle'}
	\end{equation*}
	and in particular due to the boundary conditions on $\widetilde{Q}_\ub$ we obtain
	\[C = \left.\overline{\partial_d\langle E\rangle' }\right|_{x_d=0} = - \Nu.\]
	Therefore for any $\delta \in(0,1]$ it holds
	\begin{equation*}
		- \Nu = \frac{1}{\delta} \int_0^{\delta} \overline{\partial_d \langle E \rangle' + \frac{1}{2} \left\langle \left[\widetilde{Q}_\ub \nabla E\right]^{(d)}\right\rangle'} {\rm d}x_d.
	\end{equation*}
	\noindent
	Next observe that
	\[\left\langle \left[\widetilde{Q}_\ub \nabla E\right]^{(d)}\right\rangle' = \left\langle \sum_{j=1}^d \widetilde{q}_{d,j} \partial_j E\right\rangle'.\]
	Thus we deduce, again by using Lemma \ref{lem:ExOp2}, that
    \begin{align*}
		- \Nu &= \frac{1}{\delta} \int_0^{\delta} \overline{\partial_d \langle E \rangle' + \frac{1}{2} \left\langle \left[\widetilde{Q}_\ub \nabla E\right]^{(d)}\right\rangle'} {\rm d}x_d\notag\\ 
		&\overset{\text{Lem. \ref{lem:ExOp2}}}{=} \overline{\frac{1}{\delta} \int_0^{\delta} \partial_d \langle E \rangle' + \frac{1}{2} \left\langle \left[\widetilde{Q}_\ub \nabla E\right]^{(d)}\right\rangle' {\rm d} x_d}\notag\\
		&= \overline{\frac{1}{\delta}\left(\langle E \rangle'(x_d = \delta) - \langle E \rangle'(x_d = 0)\right) + \frac{1}{2\delta} \sum_{j=1}^d \int_0^{\delta} \langle \widetilde{q}_{d,j} \partial_j E\rangle' {\rm d}x_d}\notag\\
		&\geq -\frac{1}{\delta} + \frac{1}{2\delta}\sum_{j=1}^d\overline{\int_0^{\delta} \langle \widetilde{q}_{d,j} \partial_j E\rangle' {\rm d}x_d}\notag\\
		&\geq -\frac{1}{\delta} - \frac{1}{2\delta} \sum_{j=1}^d\left(\int_0^{\delta} \langle |\widetilde{q}_{d,j}|^2 \rangle ' {\rm d}x_d \right)^{\frac{1}{2}}\overline{\left(\int_0^{\delta}\langle | \partial_j E|^2 \rangle' {\rm d}x_d\right)^{\frac{1}{2}}}\notag\\
		&\geq -\frac{1}{\delta} - \frac{1}{2\delta} \sum_{j=1}^d\left(\int_0^{\delta} \langle |\widetilde{q}_{d,j}|^2 \rangle ' {\rm d}x_d \right)^{\frac{1}{2}}\left(\overline{\int_0^{\delta}\langle | \partial_j E|^2 \rangle' {\rm d}x_d}\right)^{\frac{1}{2}}\notag\\
		&\geq -\frac{1}{\delta} - \frac{1}{2\delta} \sum_{j=1}^d\left(\int_0^{\delta} \langle |\widetilde{q}_{d,j}|^2 \rangle ' {\rm d}x_d \right)^{\frac{1}{2}}\Nu^{\frac{1}{2}}
	\end{align*}
	where we have used Jensen's inequality, the maximum principle and in particular in the last line Lemma \ref{lem:NuId} to conclude that for all $j= 1, \dots, d$ it holds
	\[\overline{\int_0^{\delta}\langle | \partial_j E|^2 \rangle' {\rm d}x_d} \leq \overline{\int_0^{\delta}\langle | \nabla  E|^2 \rangle' {\rm d}x_d} \leq \overline{\int_0^1 \left \langle \left| \left( {\rm Id} + \frac{1}{2}\widetilde{Q}_\ub\right)^{\frac{1}{2}}\nabla E\right|^2\right \rangle'{\rm d}x_d} = \Nu.\]
	\noindent
	The assertion follows by applying Young's inequality.
\end{proof}
\begin{proof}[Proof (Lemma \ref{lem:ExOp2})]
	\noindent
	Again denote $E := \mathbb{E}[T]$. Further denote
	\begin{equation*}
		f^{(t_0)}(x_d) := \frac{1}{t_0}\int_0^{t_0} \partial_d \langle E\rangle' + \left \langle \left[\widetilde{Q}_\ub \nabla E\right]^{(d)}\right \rangle' {\rm d}t
	\end{equation*}
	as well as
	\begin{equation*}
		{\rm osc}_{x_d} f := \sup_{x_d\in[0,1]} f(z) - \inf_{x_d\in[0,1]} f(z).
	\end{equation*}
	Then
	\[{\rm osc}_{x_d} f^{(t_0)} \lesssim \frac{1}{t_0}.\]
	Indeed:
	\begin{align*}
		&{\rm osc}_{x_d} f^{(t_0)} = {\rm osc}_{x_d}\left(\frac{1}{t_0} \int_0^{t_0} \partial_d \langle E \rangle' + \frac{1}{2} \left \langle \left[\widetilde{Q}_\ub \nabla E\right]^{(d)}\right \rangle' {\rm d}t \right)\\
		&= \sup_{x_d\in[0,1]}\left(\frac{1}{t_0} \int_0^{t_0} \partial_d \langle E \rangle' + \frac{1}{2} \left \langle \left[\widetilde{Q}_\ub \nabla E\right]^{(d)}\right \rangle' {\rm d}t \right) - \inf_{x_d\in[0,1]}\left(\frac{1}{t_0} \int_0^{t_0} \partial_d \langle E \rangle' + \frac{1}{2} \left \langle \left[\widetilde{Q}_\ub \nabla E\right]^{(d)}\right \rangle' {\rm d}t \right)\\
		&\leq \int_0^1 \left| \partial_d \frac{1}{t_0} \int_0^{t_0} \partial_d \langle E \rangle' + \frac{1}{2} \left \langle \left[\widetilde{Q}_\ub \nabla E\right]^{(d)}\right \rangle' {\rm d}t \right| {\rm d}x_d \\
		&\leq \sup_{x_d\in[0,1]} \left| \partial_d\frac{1}{t_0} \int_0^{t_0} \partial_d \langle E \rangle' + \frac{1}{2} \left \langle \left[\widetilde{Q}_\ub \nabla E\right]^{(d)}\right \rangle' {\rm d}t \right|\\
		&=\sup_{x_d\in[0,1]} \left| \frac{1}{t_0} \int_0^{t_0} \partial_d\left(\partial_d \langle E \rangle' + \frac{1}{2} \left \langle \left[\widetilde{Q}_\ub \nabla E\right]^{(d)}\right \rangle'\right) {\rm d}t \right|\\
		&= \sup_{x_d\in[0,1]} \left| \frac{1}{t_0} \int_0^{t_0} \partial_t \langle E \rangle' {\rm d}t\right| \lesssim \frac{1}{t_0}
	\end{align*}
	having used the evolution equation \eqref{eq:EvolExp} and the maximum principle. In particular, this implies that
	\[\limsup_{t_0 \rightarrow \infty} f^{(t_0)}(x_d) = C \quad \text{independent of }x_d.\]
	This yields pointwise convergence and existence of a dominating function, thus the second claim follows from applying the dominated convergence theorem.
   
\end{proof}
\section{Proof of Corollary \ref{cor:PeScaling}}\label{sec:boundPeclet}
Recall the representation \eqref{eq:cylinBM}
\begin{equation*}
    {\rm d} \ub = \sum_{k \in \mathbb{Z}} \lambda_k e_k(\xx) {\rm d}\beta_t^k
\end{equation*}
and that the (non-dimensional) Péclet number is defined by
\begin{equation*}
    \Pe^2 = \int_0^1 \left \langle {\rm tr}\left(\widetilde{Q}_\ub\right)\right\rangle' {\rm d}x_d
\end{equation*}
where
\[\widetilde{Q}_\ub(\xx) = \sum_{k \in \mathbb{Z}} \lambda_k^2 e_k(\xx) \otimes e_k(\xx)\]
and hence
\[{\rm tr}\left(\widetilde{Q}_\ub\right) = \sum_k \lambda_k^2 |e_k|^2.\]
Therefore we obtain
\begin{equation}\label{eq:PeLambda}
    \Pe^2 = \sum_k \lambda_k^2.
\end{equation}
Now let
\[\widetilde{Q}_0(\xx) := \sum_k \mu_k^2 e_k(\xx) \otimes e_k(\xx)\]
where
\[\mu_k := \frac{\lambda_k}{\Pe} = \frac{\lambda_k}{\sqrt{\sum_l \lambda_l^2}}.\]
Hence
\[\int_0^1 \left \langle {\rm tr}\left(\widetilde{Q}_0\right)\right \rangle' {\rm d} x_d = 1.\]

We now derive upper bounds on the vertical covariance $\tilde q_{d,j}$ which, combined with Theorem \ref{thm:NuGeneralBound}, yields the scaling of
the Nusselt number $\Nu$ in $\Pe$. Recall
 \begin{corollary}
For all $x_d\in[0,1]$,
\begin{align}
\left\langle|\widetilde q_{d,j}|^2\right\rangle'&\leq C\,\Pe^4x_d^6,
\qquad j=1,\dots,d-1, \label{eq:QDecay1}\\
\left\langle|\widetilde q_{d,d}|^2\right\rangle'&\leq C\,\Pe^4x_d^8,
\label{eq:QDecay2}
\end{align}
with $C$ depending only on $d$ and $\sum_k\mu_k^2\,\|e_k\|_{C^3}^2<\infty$. Consequently,
\[\Nu\lesssim\Pe^{2/3} \quad \mbox{ for }\quad \Pe\ge1.\]
If moreover
$\left\langle|\widetilde q_{d,j}|^2\right\rangle' = 0$ for $j<d$, then \[\Nu\lesssim\Pe^{1/2}\,,\] and this bound is sharp for $d=3$.
\end{corollary}
\begin{proof}%[Proof of Corollary \ref{cor:PeScaling}]
   
    \textbf{Upper bounds:} By regularity of the vector-valued spatial coefficients $e_k$, we may perform a Taylor expansion of the components $e_k^j, \, j=1,\dots,d,$ around $x_d = 0$
    \begin{align*}
        e_k^j(\xx',x_d) &= e_k^j(\xx',0) + x_d \left(\partial_d e_k^j\right)(\xx',0) + \frac 12 x_d^2 \left(\partial_d^2 e_k^j\right)(\xx',0) + O\left(x_d^3\right)\\
        &= x_d \left(\partial_d e_k^j\right)(\xx',0) + x_d^2 \left(\partial_d^2 e_k^j\right)(\xx',0) + O\left(x_d^3\right)
    \end{align*}
    where we have used the no-slip boundary condition. Now for the case $j= d$, together with the incompressibility condition we may further deduce
    \[\left(\partial_d e_k^d\right)(\xx',0) = - \sum_{j=1}^{d-1} \left(\partial_j e_k^j\right)(\xx',0) = 0,\]
    hence
    \[e_k^d(\xx',x_d) = \frac 12 x_d^2 \left(\partial_d^2 e_k^d\right)(\xx',0) + O(x_d^3).\]
    Recalling that
    \[\widetilde{q}_{d,j}(\xx) = \sum_k \lambda_k^2 e_k^d(\xx)e_k^j(\xx)\]
    and $\lambda_k = \Pe \mu_k$ we deduce with the above that for $j=1,\dots,d-1$ it holds
    \begin{equation*}
        \langle |\widetilde{q}_{d,j}|^2 \rangle' \lesssim \Pe^4 \left(x_d^6\left \langle \left| \sum_k \mu_k^2 \left(\partial_d^2 e_k^d\right)\left(\partial_d e_k^j\right)\right|^2\right \rangle' + O\left(x_d^8\right)\right),
    \end{equation*}
    while for $j=d$ we obtain
    \begin{equation*}
        \langle |\tilde{q}_{d,d}|^2 \rangle' \lesssim \Pe^4 \left(x_d^8\left \langle \left| \sum_k \mu_k^2 \left(\partial_d^2 e_k^d\right)^2\right|^2\right \rangle' + O\left(x_d^{10}\right)\right),
    \end{equation*}
    Here the constants in $O(\cdot)$ are uniform in $\xx'$ and depend only on
$\sum_k\mu_k^2\|e_k\|_{C^3}^2<\infty$; since $x_d^8\le x_d^6$ and $x_d^{10}\le x_d^8$ on
$[0,1]$, this proves \eqref{eq:QDecay1} and \eqref{eq:QDecay2} for all
$x_d\in[0,1]$.\\
    In view of \eqref{eq:genUpperBound}, summing the bounds
    \eqref{eq:QDecay1}-\eqref{eq:QDecay2} over $j$, we obtain for
    $\delta\in(0,1]$
    \begin{equation*}
    \Nu\lesssim\frac1\delta+\frac1{\delta^2}\int_0^\delta
    \Pe^4\left(x_d^6+x_d^8\right)\,{\rm d}x_d
    \lesssim\frac1\delta+\Pe^4\delta^5 ,
    \end{equation*}
    and choosing $\delta=\Pe^{-2/3}$, which is admissible for $\Pe\ge1$, we
    obtain $\Nu\lesssim\Pe^{2/3}$. If moreover $\langle |\widetilde{q}_{d,j}|^2 \rangle'\equiv0$ for
    $j<d$, only the term with $x_d^8$ remains, thus
    $\Nu\lesssim\delta^{-1}+\Pe^4\delta^7$, and the choice
    $\delta=\Pe^{-1/2}$ gives $\Nu\lesssim\Pe^{1/2}$.

    \smallskip    
  
      \textbf{Sharpness:}
      Here we construct a stochastic velocity field in dimension $d=3$, for which $\widetilde{q}_{i,j}=0$ unless $i=j=3$, which attains $\Nu\sim \Pe^{\frac 12}$. In this regard we define
      \begin{align*}
          a(x_3) &:= x_3^2(1-x_3)^2,\\
          e_1(\xx) &:= \frac{1}{C}\hat{e}_1(\xx), \quad \hat{e}_1(\xx) :=\begin{bmatrix}
              a'(x_3)\cos(px_1) \\ 0 \\ pa(x_3)\sin(px_1)
          \end{bmatrix},\\
          e_2(\xx) &:= \frac{1}{C}\hat{e}_2(\xx), \quad \hat{e}_2(\xx) :=\begin{bmatrix}
              a'(x_3)\sin(px_1) \\ 0 \\ -pa(x_3)\cos(px_1)
          \end{bmatrix},\\
          \lambda_1 &= \lambda_2 = \sqrt{\frac{\Pe^2}{2}}  =: \Lambda \quad \text{(cf. \eqref{eq:PeLambda})}
      \end{align*}
      where $p = \frac{2\pi}{L}$, with normalisation constant
      \[C^2= L^2 \left(\|a'\|_{L^2([0,1])}^2 + p^2\|a\|_{L^2([0,1])}^2\right).\]
      One can easily check that $e_k$ are $L^2(\D)$-orthonormal, divergence-free and satisfy the boundary conditions $e_k(x_3=0)=e_k(x_3=1)=0$. Thus for two one-dimensional independent Brownian motions $\beta^1, \beta^2$ we consider
      \begin{equation}\label{eq:uSharp}
          \ub = \sum_{k=1}^2 \lambda_k e_k \beta^k
      \end{equation}
      with covariance kernel
      \[Q(\xx,\yy) = \sum_{k=1}^2 \lambda_k^2 e_k(\xx) \otimes e_k(\yy) = \frac{\Lambda^2}{C^2}\sum_{k=1}^2 \hat{e}_k(\xx) \otimes \hat{e}_k(\yy).\]
      With the above choice we find that
      \begin{align*}
          \widetilde{q}_{3,1}(\xx) &= \frac{\Lambda^2}{C^2}\left(p a(x_3) a'(x_3) \cos(px_1)\sin(x_1) - p a(x_3) a'(x_3)\cos(px_1)\sin(px_1)\right) = 0,\\
          \widetilde{q}_{3,2}(\xx) &= 0,\\
          \widetilde{q}_{3,3}(\xx) &= \frac{\Lambda^2}{C^2}(pa(x_3))^2.
      \end{align*}
      Hence $\widetilde{q}_{3,j}$ only depend on $x_3$, thus as discussed before Lemma \ref{lem:NuReprZdep}, we may therefore use its result \eqref{eq:NuIdentity}, namely that in this case it holds 
      \[\Nu = \left(\int_0^1 \frac{1}{1+ \frac{1}{2}\widetilde{q}_{3,3}(x)}{\rm d} x\right)^{-1}.\]
      Thus we obtain
      \[\Nu = \left(\int_0^1 \frac{1}{1+ \frac{1}{2}\frac{\Lambda^2}{C^2}\left(pa(x)\right)^2}{\rm d} x\right)^{-1} = \left(\int_0^1 \frac{1}{1+ \tilde{C}\Pe^2x^4(1-x)^4}{\rm d} x\right)^{-1},\]
      where $\tilde{C}=\frac{p^2}{4C^2}$. Now observe that it holds
  \[\int_0^1 \frac{1}{1+ \tilde{C}\Pe^2x^4(1-x)^4}{\rm d} x = 2 \int_0^\frac{1}{2} \frac{1}{1+ \tilde{C}\Pe^2x^4(1-x)^4}{\rm d} x.\]
  Next we estimate
  \begin{align*}
      \int_0^\frac{1}{2} \frac{1}{1+ \tilde{C}\Pe^2x^4(1-x)^4}{\rm d} x &\leq \int_0^\frac{1}{2} \frac{1}{1+ \frac{\tilde{C}}{16}\Pe^2x^4}{\rm d} x \\
      &\leq \int_0^\infty \frac{1}{1+ \frac{\tilde{C}}{16}\Pe^2x^4}{\rm d} x \\
      &= \frac{2}{\tilde{C}^{\frac{1}{4}}\Pe^\frac{1}{2}}\int_0^\infty \frac{z^{-\frac{3}{4}}}{1+ z}{\rm d} z = \frac{2\sqrt{2}\pi}{\tilde{C}^{\frac{1}{4}}} \Pe^{-\frac{1}{2}}.
  \end{align*}
  This provides the lower bound
  \[\Nu \gtrsim \Pe^{\frac{1}{2}}.\]
  Let us further note that in this example we may even determine the upper bound: it holds
  \begin{align*}
   \int_0^\frac{1}{2} \frac{1}{1+ \tilde{C}\Pe^2x^4(1-x)^4}{\rm d} x&\geq \int_0^{\frac 12}\frac{1}{1+\tilde{C}\Pe^2 x^4} {\rm d}x\\
   &= 
   \frac{1}{\tilde{C}^{\frac{1}{4}}\Pe^{\frac 12}}\int_{0}^{\frac{1}{2}\tilde{C}^{\frac{1}{4}}\Pe^{\frac 12}}\frac{1}{1+ z^4} {\rm d}z\\
   &= \frac{1}{\tilde{C}^{\frac{1}{4}}\Pe^{\frac 12}}\left[\int_{0}^{\infty}\frac{1}{1+ z^4} {\rm d}z- \int_{\frac{1}{2}\tilde{C}^{\frac 14}\Pe^{\frac 12}}^{\infty}\frac{1}{1+ z^4}{\rm d}z\right]\\
   &\geq\frac{1}{\tilde{C}^{\frac{1}{4}}\Pe^{\frac 12}}\left[\frac{\sqrt{2}\pi}{4}-\frac{1}{3\tilde{C}^{\frac 34}\Pe^{\frac 32}}\right]\,,
   \end{align*}
   and the last bracket is bigger than $\frac{\sqrt{2}\pi}{8}$ as soon as $\Pe\geq \left(\frac{8}{3\sqrt{2}\pi}\right)^{\frac 23}\tilde{C}^{-\frac{1}{2} }$, yielding
   \[ \int_0^{\frac 12}\frac{1}{1+\tilde{C}\Pe^2 x^4} {\rm d}x \geq \frac{\sqrt{2}\pi}{8}\frac{1}{\tilde{C}^{\frac{1}{4}}\Pe^{\frac 12}}\]
   and hence the upper bound
  \[\Nu \lesssim \Pe^{\frac{1}{2}}\]
  where the constant can be made independent of $L$ by using that
  \begin{align*}
      \tilde{C} = \frac{p^2}{4C^2} &= \frac{p^2}{4L^2\left(\|a'\|_{L^2([0,1])}^2 + p^2 \|a\|_{L^2([0,1])}^2\right)} \leq \frac{1}{4\|a\|_{L^2([0,1])}^2}
  \end{align*}
  for $L\geq 1$.
\end{proof}

\appendix
\section{Appendix}\label{app:rescaling}
In this section we use the stochastic framework introduced in Section \ref{sec:stochFramework}: in the dimensional setting, let $\Theta$ be a probabilistically strong, mild solution to 
\begin{equation*}
{\rm d}\Theta+\left(\circ\sqrt{\mathcal{Q}_{\vv}}\,{\rm d}W\right)\cdot\nabla\Theta
=\kappa\Delta\Theta\,{\rm d}t
\end{equation*}
where
\[\circ\sqrt{\mathcal{Q}_{\vv}}\,{\rm d}W = \sum_k \gamma_k f_k \circ {\rm d}B^k.\]
That is, for every $t \geq 0$ and almost every $\xx$ we may write
\begin{align*}
    \Theta(t,\xx) = \Theta(0,\xx) &- \sum_k \int_0^t \gamma_k (f_k(\xx) \cdot \nabla )\Theta(s,\xx) \circ {\rm d}B_s^k\\
    &+ \kappa \int_0^t \Delta \Theta (s,\xx) {\rm d}s.
\end{align*}
We introduce the following change of variables $\xx = h \hat{\xx}, t = \frac{h^2}{\kappa}\hat{t}$ and introduce $\hat{\Theta}(\hat{t},\hat{\xx}) = \Theta(\frac{h^2}{\kappa}\hat{t}, h\hat{\xx})$. This change of variables induces $\nabla = h^{-1}\hat{\nabla}, \Delta = h^{-2} \hat{\Delta}$; furthermore due to the scale invariance property of Brownian motion we obtain that it holds in law
\begin{equation}\label{eq:BMscaleInv}
    B_t^k \overset{d}{=} \frac{h}{\sqrt{\kappa}} B^k_{\hat{t}}.
\end{equation}
Thus in distribution, $\hat{\Theta}$ is equivalent to a solution of (with abuse of notation)
\begin{equation*}
{\rm d}\hat{\Theta}+\left(\circ\sqrt{\mathcal{Q}_{\ub}}\,{\rm d}\hat{W}\right)\cdot\hat{\nabla}\hat{\Theta}
=\hat{\Delta}\hat{\Theta}\,{\rm d}\hat{t}
\end{equation*}
where
\[\circ\sqrt{\mathcal{Q}_{\ub}}\,{\rm d}\hat{W} = \sum_k \frac{\gamma_k}{\sqrt{\kappa}} f_k(h \cdot)\circ {\rm d}B^k =: \sum_k \lambda_k e_k \circ {\rm d}\beta^k.\]
Therefore we obtain the following relation between the covariance kernels
\begin{equation*}
Q_{\ub}(\hat{\boldsymbol x},\hat{\boldsymbol y})
=\frac1\kappa\,Q_{\vv}(h\hat{\boldsymbol x},h\hat{\boldsymbol y}),
\qquad\text{in particular}\qquad
\widetilde Q_{\ub}(\hat{\boldsymbol x})=\frac1\kappa\,\widetilde Q_{\vv}(h\hat{\boldsymbol x}).
\end{equation*}
Note that the noise does not rescale like a velocity. For a
deterministic velocity, the full time factor enters, ${\rm d}\boldsymbol
x/{\rm d}t=(\kappa/h)\,{\rm d}\hat{\boldsymbol x}/{\rm d}\hat t$, so that
$\ub=h\vv/\kappa$. For the noise, only its square-root enters through
\eqref{eq:BMscaleInv}, such that the covariance
acquires the factor $1/\kappa$ rather than $h^2/\kappa^2$. In
dimensionless variables the covariance budget therefore reads
\begin{equation*}
\Pe^2=\int_0^1\left\langle\mathrm{tr}\left(\widetilde Q_{\ub}\right)\right\rangle'\,{\rm d}\hat{x}_d .
\end{equation*}

\section*{AI Usage Disclosure}
All results in this work, including their proofs, were obtained by the authors without the use of any AI tool. During the preparation of the manuscript, the authors used Claude (Anthropic) to improve the language, style and presentation of the introduction and abstract, and to assist in interpreting the results of \cite{HCD2014}. All suggestions were critically assessed by the authors, who take full responsibility for the content of the publication.

\section*{Acknowledgements}
During this project, FF and TL have received funding from the European Research Council (ERC) under the EU-Horizon Europe ERC-2021-ADG research and innovation programme (project ``Noise in Fluids", grant agreement no. 101053472), while TL has received further funding from the Deutsche Forschungsgemeinschaft (DFG, German Research Foundation) - project number 569224218. CN gratefully acknowledges the support of TRR181 and the Scuola Normale Superiore (Pisa), which enabled mutual visits in 2021–2022, during which the project was initiated in collaboration with FF.

\bibliographystyle{plain}
\bibliography{bibliography}

\end{document}